\documentstyle[11pt]{article}
   \csname @addtoreset\endcsname{equation}{section}
    \csname @addtoreset\endcsname{figure}{section}
    \csname @addtoreset\endcsname{table}{section}

\newcounter{thanksnum}
\def\thanksnumber#1
{\setcounter{thanksnum}{\value{footnote}}\setcounter{footnote}{#1}%
                     \addtocounter{footnote}{-1}\footnotemark
                     \setcounter{footnote}{\value{thanksnum}}}
\def\newtheoremz#1{\@ifnextchar[{\@othmz{#1}}{\@nthmz{#1}}}

\def\@nthmz#1#2{%
\@ifnextchar[{\@xnthmz{#1}{#2}}{\@ynthmz{#1}{#2}}}

\def\@xnthmz#1#2[#3]{\expandafter\@ifdefinable\csname #1\endcsname
{\@definecounter{#1}\@addtoreset{#1}{#3}%
\expandafter\xdef\csname the#1\endcsname{\expandafter\noexpand
  \csname the#3\endcsname \@thmcountersepz \@thmcounterz{#1}}%
\global\@namedef{#1}{\@thmz{#1}{#2}}\global\@namedef{end#1}{\@endtheoremz}}}

\def\@ynthmz#1#2{\expandafter\@ifdefinable\csname #1\endcsname
{\@definecounter{#1}%
\expandafter\xdef\csname the#1\endcsname{\@thmcounterz{#1}}%
\global\@namedef{#1}{\@thm{#1}{#2}}\global\@namedef{end#1}{\@endtheoremz}}}

\def\@othmz#1[#2]#3{\expandafter\@ifdefinable\csname #1\endcsname
  {\global\@namedef{the#1}{\@nameuse{the#2}}%
\global\@namedef{#1}{\@thmz{#2}{#3}}%
\global\@namedef{end#1}{\@endtheoremz}}}

\def\@thmz#1#2{\refstepcounter
    {#1}\@ifnextchar[{\@ythmz{#1}{#2}}{\@xthmz{#1}{#2}}}

\def\@xthmz#1#2{\@begintheoremz{#2}{\csname the#1\endcsname}\ignorespaces}
\def\@ythmz#1#2[#3]{\@opargbegintheoremz{#2}{\csname
       the#1\endcsname}{#3}\ignorespaces}

\def\@thmcounterz#1{\noexpand\arabic{#1}}
\def\@thmcountersepz{.}
\def\@begintheoremz#1#2{ \trivlist \item[\hskip \labelsep{\bf #1\ #2}]}
\def\@opargbegintheoremz#1#2#3{ \trivlist
      \item[\hskip \labelsep{\bf #1\ #2\ (#3)}]}
\def\@endtheoremz{\endtrivlist}

\newtheorem{theorem}{Theorem}[section]
\newtheorem{lemma}{Lemma}[section]

\newtheorem{proposition}{Proposition}[section]

\newtheorem{definition}{Definition}[section]

\def\e{\varepsilon}

\def\defi{\stackrel{{\scriptscriptstyle \Delta}}{=}}

\def\a{\alpha}
\def\d{\delta}
\def\o{\omega}
\def\O{\Omega}

\def\Y{{\cal Y}}
\def\F{{\cal F}}
\def\w{\widehat}

\def\esssup{\mathop{\rm ess\, sup}}
\def\const{{\rm const\,}}

\def\R{{\bf R}}
\def\E{{\bf E}}
\def\P{{\bf P}}

\def\H{{\cal H}}

\def\L{L}

\def\b{\beta}
\def\s{\delta}
\def\g{\gamma}

\def\ww{\widetilde}
\def\X{{\cal X}}

\def\oo{\bar}
\def\s{\sigma}

\def\p{\partial}

\def\L{{\cal L}}

\newcommand{\be}{\begin{equation}}
\newcommand{\ee}{\end{equation}}
\newcommand{\bd}{\begin{displaymath}}
\newcommand{\ed}{\end{displaymath}}
\newcommand{\ba}{\begin{array}{ll}}
\newcommand{\ea}{\end{array}}
\font\sm=cmr10

\def\BS{{\scriptscriptstyle BS}}

\def\Re{{\rm Re\,}}

\def\ww{\tilde}

\def\H{{\cal H}}
\def\u{ u }
\def\v{{\rm\bf v }}
\title{Pricing rule
 based on non-arbitrage arguments for random volatility and volatility smile}
\author{N.G.Dokuchaev\thanks{Corresponding address:
 Department of Mathematics and Computer Science, The University
of West Indies, Mona, Kingston 7, Jamaica W.I.
 Email: ndokuch@uwimona.edu.jm}
\\
 {\sm   Department of Mathematics and Computer Science,} {\sm The
University of West Indies, Mona}\\{\sm  Jamaica}}
\begin{document}
\maketitle
\begin{abstract}
We consider a generic market model with a single stock and with
random volatility. We assume that there is a number of tradable
options for that stock with different strike prices. The paper
states the problem of finding  a pricing rule that gives
Black-Scholes price for at-money options and such that the market
is arbitrage free for any number of tradable options, even if
there are two Brownian motions only: one drives the stock price,
the other drives the volatility process. This problem is reduced
to solving a parabolic equation.
\end{abstract}
\section{Introduction}
We consider a market model with a single stock and with random
volatility. Most practitioners have adapted the famous
Black-Scholes model as the premier model for pricing and hedging
of options. This model consists of two assets: the risk free bond
or bank account and the risky stock. It is assumed that the
dynamics of the stock is given by a random process with some
standard deviation of the stock returns (the volatility
coefficient, or volatility). The dynamics of bonds is
deterministic and exponentially increasing with a given risk-free
rate. In the classic  Black-Scholes model, the volatility is
assumed to be given and fixed. However, empirical research shows
that the real volatility is time-varying and random. Moreover, it
is commonly recognized  that Black-Scholes formula gives unbiased
estimation for at-money options only, and it gives a systematic
error for in-money and out-of-money options; in fact, that means
that  there is a gap between historical and implied volatility
that  generates so-called {\it volatility smile} for the implied
volatility; see e.g. Black and Scholes (1972), Day and Levis
(1992), Derman {\it et al.} (1996), Hauser and Lauterbach (1997),
Taylor and  Xu (1994). A very detailed review
 can be found in Mayhew (1995)). Many authors emphasize that the
main difficulty in modifying the Black--Scholes and Merton models
is taking into account  this fact.
\par
To fill this gap, a number of deterministic and stochastic equations for volatility
were proposed (see e.g. Christie (1982), Johnson  and  Shanno (1987), Hull  and  White
(1987), Masi {\it et al.} (1994) and more recent papers in Jarrow (ed.) (1998). In
some other approach, a special temporal scale is used to find the time when historical
volatility coincides with implied volatility (see e.g. Geman  and Ane (1996)).
\par
This paper suggests some another approach to the problem based on non-arbitrage arguments
(see, e.g., Harrison and Pliska (1981) and  Jouini (1996)). We take
a very generic model for the random volatility process (more
precisely, for historical volatility process), and then look for a
pricing rule that  \begin{itemize}\item gives Black-Scholes price
for at-money options;
\item ensures  that  arbitrage possibilities are absent for an extended market
that includes a
number of tradable options for that stock with different strike
prices, and such that there are two Brownian motions only: one drives the
stock price, the other drives the volatility process.
\end{itemize}
Notice that similar non-arbitrage arguments are used traditionally
for bond pricing; for the generic model, there is only one Brownian motions driving
the interest rate there are many bonds. We explore the same approach for multi-option market to find
 a model for the volatility smile.
\par
We show below the required pricing rule exists if
some boundary value problem for a parabolic equation is solvable; the price is expressed via
solution of this equation. This boundary value problem  looks simple,
but, unfortunately, it appears to be difficult to solve; the boundary
conditions cannot be reduced to known ones. Most part
of the paper devoted analysis of solvability of this boundary value problem. We
obtained some existence and uniqueness theorem and prior estimations for the solution;
we found a also formula
for the Fourier transform of the solution via expectation of a
integral functional of the volatility  process.
\section{Definitions} Consider the diffusion model of a
securities market  consisting of a risk free bond or bank account
with the price $B(t), $ ${t\ge 0}$, and
 a risky stock with price $S(t)$, ${t\ge 0}$. The prices of the stocks evolve
 according to the following
stochastic differential equation \be \label{S}
dS(t)=S(t)\left(a(t)dt+\s(t) dw(t)\right), \quad t>0, \ee where
$w(t)$ is a standard Wiener process,  $a(t)$ is an appreciation
rate, $\s(t)$ is a volatility coefficient.  The initial price
$S(0)>0$ is a  given deterministic constant.  The price of the
bond evolves according to the following equation
\begin{equation}
\label{B}
B(t)=e^{rt}B(0),
\end{equation}
where $r\ge 0$ and $B(0)$ are  given constants. \par Set
$v(t)\defi \s(t)^2$. We assume that the process for $v(t)$ evolves
as \be \label{v} dv(t)=v(t)\left(\w a dt+\w\s d\w w(t)\right),
\quad t>0, \ee where $\w a$ and $\w \s\neq 0$ are known constants,
and $\w w(\cdot)$ is a Wiener process independent of $w(\cdot)$.
The initial volatility $v(0)>0$ is a  given deterministic
constant. Moreover, we assume for the sake of simplicity that the
volatility process has zero drift, i.e. $\w a=0$ (an extension for
the case $\w a\neq 0$ does requires principal changes and is
rather technical).
\par
We  assume that $(w(\cdot),\w  w(\cdot))$ is a standard Wiener
process on  a given standard probability space $(\O,\F,\P)$, where
$\O$ is a set of elementary events, $\F$ is a  complete
$\s$-algebra of events, and $\P$ is a probability measure.
\par
Further, we assume that there are available European options on
that stocks with the same expiration time $T$ and different strike
prices $K\in {\cal K}$, where ${\cal K}\subset(0,+\infty)$ is a
given set. For simplicity, we assume that the set ${\cal K}$ is
finite, i.e. ${\cal K}=\{K_1,...,K_N\}$, where $N$ is an integer,
possibly a large number. (Generalization for the case of infinite
number of $K$ and different expiration times is rather
technical).\par We shall denote $P_c(t,K)$ and $P_p(t,K)$ the
prices for the call option with the claim $(S(T)-K)^+$ and for the
put option with the claim $(K-S(T))^+$ respectively, where $T$ is
the expiration time, $K$ is the strike price.  We shall consider
options as additional tradable assets, i.e. we shall consider {\it
bond-stock-options market}. Under our assumptions, the market is
incomplete (i.e. options can not be replicated). We shall look for
a reasonable model for prices $\{P_c(t,K),P_p(t,K)\}_{K\in{\cal
K}}$ (i.e. a reasonable pricing rule); the main requirement is
that the market must be arbitrage free. Since we have only two
driving Brownian motions $w(\cdot)$ and $\w w(\cdot)$, we need to
prevent arbitrage possibilities for the case of  large $N$.
\par
Let $H_{\BS,c}(x,v,t,K)$ and $H_{\BS,p}(x,v,t,K)$ denotes
Black-Scholes prices for the put and call options with the claims
$(S(T)-K)^+$ and $(K-S(T))^+$ respectively given condition
$S(t)=x$, $\s(s)^2|_{[t,T]}\equiv v$. (Note that
$\E\{\s(s)^2|v(t)=v\}$ is $v$ for any $s\ge t$.) More precisely,
\be\label{BS}\ba
 \ww H_{\BS,c}(x,v,t,\ww
K)=x\Phi(d_+(x,v,t,\ww K ))-\ww K\Phi(d_+(x,v,t,\ww K )), \\
\ww H_{\BS,p}(x,v,t,\ww K)= \ww H_{\BS,1}(x,v,t,\ww K)+x-\ww K, \ea \ee where
\be\label{dpm}\ba
 \Phi(x)
\defi\frac{1}{\sqrt{2\pi}}\int_{-\infty}^x
e^{-\frac{t^2}{2}}dt, \\ d_+(x,v,t,\ww K )\defi\frac{\log{x}-\log{\ww
K}}{\sqrt{(T-t)v}}+ \frac{\sqrt{(T-t)v}}{2}, \quad d_-(x,v,t,\ww K
)\defi\frac{\log{x}-\log{\ww K}}{\sqrt{(T-t)v}}- \frac{\sqrt{(T-t)v}}{2}. \ea \ee
\par
To make the following definitions more compact, we shall denote $$ \ba
P_1(t,K)\defi P_c(t,K), \quad P_2(t,K)\defi
P_p(t,K),\\
H_{\BS,1}(x,v,t,K)\defi H_{\BS,c}(x,v,t,K),\quad
H_{\BS,2}(x,v,t,K)\defi H_{\BS,p}(x,v,t,K).\ea$$
\begin{definition} A process $v_K(t)$ is said to be {\it implied volatility}   given the strike price
$K$ and the type of option (i.e. put or call), if $P_j(t,K)\equiv
H_{\BS,j}(t,S(t),v_{K,j}(t),K)$, $j=1,2$.
\end{definition}
Clearly, the problem of finding prices
$\{P_c(t,K),P_p(t,K)\}_{K\in{\cal K}}$  can be reformulated as the
problem of finding implied volatilities $\{v_{K,j}(t)\}_{K\in{\cal
K},i=1,2}$.
\par
As known, if $v$ is non-random,  i.e. if $\w\s=0$, then the market
with Black-Scholes prices $P_j(t,K)$$\defi
H_{\BS,j}(t,S(t),v(t),K)$ is arbitrage free (i.e. the market is
arbitrage free with $v_K(t)\equiv v$, when the implied volatility
and the historical volatility are equal). Unfortunately,  for
random volatilites, this pricing rule leads to a model with
arbitrage possibilities.
\begin{proposition}\label{propBS}
If $\w\s\neq 0$ then the pricing rule   $P_j(t,K)\defi
H_{\BS,j}(t,S(t),v(t),K)$ does not allow arbitrge for $N=1$ and
allows arbitrage for $N\ge 2$.
\end{proposition}
 One of possible ways to construct an arbitrage free market
is to include  additional Brownian motions  for the equations for
the implied volatilities $v_{K,j}(t)$. (Remind  that the system
$\{v_{K,j}(t)\}$ defines prices as $P_j(t,K)=
H_{\BS,j}(t,S(t),v_K(t),K)$). For example, the model will be
arbitrage free if
  $$ dv_{K,j}(t)=v_{K,j}(t)\left(\w a_{K,j} dt+
  \w\s_{K,j} d\w w(t)+\s_{K,j}dw_{K,j}(t)\right),
\quad t>0,\quad j=1,2, $$ where $\{w_{K,j}(\cdot)\}_{K\in{\cal
K},j=1,2}$ is a system of independent Brownian motions that does
not depend on $(w(\cdot),\w w(\cdot))$, and where
$\w\s(K,j),\s_{K,j}\neq 0$ are constants.
\par
We introduce here another arbitrage free market model with only
two driving Brownian motions and arbitrarily large $N>0$.
\par
Set
$$
\ww
S(t)\defi e^{-rt}S(t), \quad \ww
K\defi e^{-rT}K,\quad\ww
P(t,\ww K)\defi e^{-rt}P(t,K), \quad K\in{\cal K}.
$$
\par
Let $\ww{\cal K}\defi\{\ww K=e^{-rT}K:\ K\in{\cal K}\}$. Let the
function $\ww H_\BS(\cdot):Q\times\ww{\cal K}\to\R$ be such that
$e^{rt}\ww H_{\BS,j}(x,v,t,\ww K)\equiv H_{\BS,j}(e^{rT}x,v,t,K)$.
\par
 Let $\P_*$ be the risk-neutral measure such that the process
 $w_*(t)\defi \int_0^t \s(s)^{-1}(a(s)-r)ds +w(t)$
 is condionally a Wiener process given $\s(\cdot)$. Let $\E_*$ be
 the corresponding expectation.
 As is known, $\ww S(t)$ is a martingale under $\P_*$, and
  $\E_*\{S(T)\,|\,S(t)\}=e^{r(T-t)}S(t)$. Thus, it is
natural to say that an option with the strike price $K$ is
at-money at time $t$ if $e^{r(T-t)}S(t)=K$, i.e.
$e^{r(T-t)}e^{rt}\ww S(t)=K$, or $\ww S(t)=\ww K$.
\par Let $\F_t$
be the filtration generates by $(S(t),v(t))$.
\par
Set $\w P_j(t,x,\ww K)\defi\E\left\{\ww P_j(t,\ww K)\,\bigl|\,\ww
S(t)=x\right\}$.
\par
We are looking for the set of prices $\{P_j(t,K)\}_{K\in{\cal
K},j=1,2}$ or $\{\ww P_j(t,\ww K)\}_{\ww K\in{\cal K},j=1,2}$ such
that the following conditions are satisfied:
\begin{itemize}
\item[(A1)] Processes $P_j(t,K)$ and  $v_{K,j}(t)$ are $\F_t$-adapted , i.e.
the only source of randomness in the model is
the pair of Brownian motions  $(w(\cdot),\w w(\cdot))$;
\item[(A2)] The market is arbitrage free
(even if the number of tradable options $N$ is arbitrarily large);
\item[(A3)] For at-money options, when $\ww S(t)\sim \ww K$, the price
is close to the Black-Scholes price; more precisely,
$$\ba \w
P_j(t,\ww K,\ww K)\equiv\ww H_{\BS,j}(\ww K,v(t),t,\ww K), \\
\left.\frac{\p\w P_j}{\p x }(t,x,\ww K)\right|_{x=\ww
K}\equiv\left.\frac{\p \ww H_{\BS,j}}{\p x}(x,v(t),t,\ww
K)\right|_{x=\ww K}. \ea
$$
In particular, if $\ww S(t)=\ww K$, then  the implied volatility
and the historical volatility are equal, i.e. $v_{K,j}(t)\equiv
v(t)$.
\item[(A4)] There exists a sequence of Markov times $\{T_k\}_{k=1}^{+\infty}$ with respect for the
filtration $\F_t$  such that $T_k\to T$ a.s.,
$\E\int_0^{T_k}P_j(t,K_i)^2dt<+\infty$, and there exists Ito's
differential $d_t P_j(t,K)$, $t\in(0,T_k)$,  for all $k$.
\end{itemize}
\subsection*{Strategies for bond-stock-options market}
Let $X(0)>0$ be the initial wealth at time $t=0$
and let $X(t)$ be the wealth at time $t>0$.
\\
We assume that
the wealth $X(t)$ at time $t\ge 0$ is
\begin{equation}
\label{X}
X(t)=\b(t)B(t)+\g_0(t)S(t)+\sum_{j=1}^2\sum_{i=1}^N\g_{i,j}(t)P_j(t,K_i).
\end{equation}
Here $\b(t)$ is the quantity of the bond portfolio, $\g_0(t)$ is the quantity of the
stock  portfolio, $\g_{i,j}(t)$ is the quantity of the options portfolio with the
strike price $K_i$ (put for $j=1$ and call for $j=2$), and
$\g(t)=\left(\g_0(t),\g_{1,1}(t),\ldots ,\g_{N,1}(t),\g_{1,2}(t),\ldots ,\g_{N,2}(t)
\right)$, $t\ge 0$. The pair $(\b(t), \g(t))$ describes the state of the bond-stocks
securities portfolio at time $t$. Each of  these pairs is  called a strategy.
\par
The process $ \ww X(t)\defi e^{-rt}X(t)$ is said to be the normalized wealth.
\begin{definition}
\label{adm} A pair $(\b(\cdot),\g(\cdot))$  is said to be an
admissible strategy if $\b(t)$, $\g_{i,j}(t)$, $j=1,2$,
$i=0,1,\ldots ,{N}$, are  random processes which are progressively
measurable with respect to the filtration $\F_t$ and such that
there exists a sequence of Markov times $\{T_k\}_{k=1}^{+\infty}$
with respect for the filtration $\F_t$  such that $T_k\to T$ a.s.
and
$$
\E\int_0^{T_k}\left(\b(t)^2dt+
S(t)^2\g_0(t)^2+\sum_{j=1}^2\sum_{i=1}^N
P_j(t,K_i)^2\g_i(t)^2\right)dt<+\infty\qquad \forall k.
$$
\end{definition}
Note that more simple static options strategies were considered in Dokuchaev (2002),
Ch.3.
\begin{definition}
A pair $(\b(\cdot),\g(\cdot))$  is said to be an admissible self-financing
strategy, if
\begin{equation}
\label{self}
dX(t)=\b(t)dB(t)+\g_0(t)dS(t)+\sum_{j=1}^2\sum_{i=1}^N\g_id_t
P_j(t,K_i).
\end{equation}
\end{definition}
In fact, (\ref{self}) is equivalent to
$$
d\ww X(t)=\g_0(t)d\ww S(t)+\sum_{j=1}^2\sum_{i=1}^N\g_i(t)d_t\ww
P_j(t,\ww K_i).
$$
(It can be shown, for example, similarly Dokuchaev and Zhou
(2001), where the same equation was derived bond-stock market).
 \par
  Let ${\bf D}\subseteq \R^m$ be
 any domain.
 Let ${\bf r}(\cdot):{\bf D}\to\R$ be a measurable function such that $\oo r(x)\ge 0$ ($\forall x$).
 Consider the Hilbert
space $L_2({\bf D},{\bf r})$ with the weight ${\bf r}(\cdot)$ that
consist of functions $f:{\bf D}\to\R$ such that $\|f(\cdot)\|_{L_2({\bf
D},{\bf r})}\defi \left(\int_{{\bf D}}|f(x)|^2{\bf r}
(x)dx\right)^{1/2}<+\infty$. Consider the Sobolev weighted spaces
$W_2^k({\bf D},{\bf r})$ with the weight ${\bf r}(\cdot)$ that
consist of functions $f(\cdot)$ such that
$\sum_{l=0}^k\|\frac{d^lf}{dx^l}(\cdot) \|_{L_2({\bf D},{\bf
r})}<+\infty$.
\par
 Let $\w r(z,y)\defi
\frac{1}{1+e^{3y}}$. Let $\H\defi L_2(\R^2,\w r)$, $\H^1\defi
W_2^1(\R^2,\w r)$.  Consider the Banach space $\Y$ of functions
$\psi(\cdot):\R\times\R\times [0,T]\to\R$ such that
\begin{itemize}\item
  $\psi(z,y,t)$ is continuous in
$(y,t)$ given $z$; \item the derivative $\frac{\p \psi}{\p
z}(z,y,t)$ is continuous in $z$ for a.e. $(y,t)$;
\item  the following estimate is satisfied:
$$
\ba \|\psi(\cdot)\|_\Y\defi
\esssup_{y,t}\left(|\psi(z,y,t)|+\left|\frac{\p \psi}{\p
z}(z,y,t)\right|\right)\\ \hphantom{xxxxxxxxxxxxxxxxxxxx}+ \sup_t
\|\psi(\cdot,t)\|_{\cal H}
+\left(\int_0^T\|\psi(\cdot,t)\|^2_{\H^1}\right)^{1/2}<+\infty.
\ea$$
\end{itemize} In particular, if $\psi(\cdot)\in\Y$ then $\frac{\p
\psi}{\p z}(z,y,t)\in L_2([0,T],\H)$ and $\frac{\p \psi}{\p
y}(z,y,t)\in L_2([0,T],\H)$.
\par
 Let $D\defi(0,+\infty)\times(0,+\infty)$, $Q\defi D\times[0,T]$.
\par  For $k>0$, let $\Y_1(k)$ denotes the class of all
functions $G=G(x,v,t):Q\to\R$ such that
   $G(x,v,t)=\psi(k\log x,\log
  y,t)$, where $\psi(\cdot)\in \Y$.
  \par Let $\Y_0(k)$ denotes the class of all functions
  $H=H(x,v,t):Q\to\R$ such that
  there exists a function $G(\cdot)\in \Y_1(k)$ such that $H(x,v,t)=G(x,v,t)+H_\BS(x,v,t,k)$.
\par Let $\varphi_1(x,k)\defi(k-x)^+$, and let
$\varphi_2(x,\ww k)\defi(x-k)^+$.
\begin{theorem}
\label{ThNAR} Assume that for any $K\in {\cal K}$, $\ww
K=e^{-rT}K$, $\varphi\equiv\varphi_j$, $j=1,2$, there exists a
function  $H(\cdot,\ww K)=H_j(\cdot,\ww K)\in\Y_0(\ww K)$ that
satisfies \be \label{parab1} \left\{ \ba \frac{\p H}{\p
t}(x,v,t,\ww K) +\frac{1}{2} x^2v\frac{\p^2 H}{\p x^2}(x,v,t,\ww
K)+\frac{1}{2} \w\s^2v^2\frac{\p^2 H}{\p v^2}(x,v,t,\ww K)=0,
\\
H(x,v,T,\ww K)=\varphi(x,\ww K),\\
H(\ww K,v,t,\ww K)=\ww H_\BS(\ww K,v,t,\ww K),\\
\left.\frac{\p  H}{\p x}(x,v,t,\ww K)\right|_{x=\w
K}=\left.\frac{\p \ww H_\BS}{\p x}(x,v,t,\ww K)\right|_{x=\w
K}.\ea \right. \ee Then conditions (A1)-(A3) are satisfied for the
market model with $\ww P_j(t,\ww K)\defi H(\ww S(t),v(t),t,\ww K)
=H_j(\ww S(t),v(t),t,\ww K)$ (in particular, the market is
arbitrage free).
\end{theorem}
{\it Remark.} We mean that (\ref{parab1}) is satisfied for
$H(\cdot,\ww K)\in\Y_0(\ww K)$ in generalized sense as equality of
generalized functions; the properties of $\Y_0(\ww K)$ are
sufficient for  correctness of all equations in (\ref{parab1}).
\par
Let $K\in {\cal K}$ be given, $\ww K=e^{-rT}K$. Set \be
\label{varphi}\ba
\varphi_0(x,v,t)=\frac{x}{\sqrt{2\pi}}e^{-\frac{d_+(x,v,t,\ww K
)^2}{2}}\biggl(\left[-\frac{\log{x}-\log{\ww
K}}{2v\sqrt{(T-t)v}}+\frac{\sqrt{(T-t)}}{4\sqrt{v}}\right]^2
+\frac{3}{2}\frac{\log{x}-\log{\ww
K}}{2v^2\sqrt{(T-t)v}}-\frac{\sqrt{(T-t)}}{8v\sqrt{v}}
\biggr)\\
\hphantom{xxx}-\frac{\ww K}{\sqrt{2\pi}}e^{-\frac{ d_-(x,v,t,\ww K
)^2}{2}}\biggl(\left[-\frac{\log{x}-\log{\ww
K}}{2v\sqrt{(T-t)v}}-\frac{\sqrt{(T-t)}}{4\sqrt{v}}\right]^2+
\frac{3}{2}\frac{\log{x}-\log{\ww
K}}{2v^2\sqrt{(T-t)v}}+\frac{\sqrt{(T-t)}}{8v\sqrt{v}}\biggl),\\
\varphi(x,v,t)\defi\frac{1}{2}\w\s^2v^2\varphi_0(x,v,t), \ea \ee
where $d_+(\cdot)$ and $d_-(\cdot)$ are  defined by (\ref{dpm}).
\par
Set
 $$
{\bf f}(z,v,t)\defi \varphi(\ww Ke^z,v,t),
$$
Set $I_1\defi(-\infty,0)$, $I_2\defi (0,+\infty)$,
\par\be\label{Ff} {\bf F}_m(\o,v,t)\defi\frac{1}{\sqrt{2\pi}}\int_{I_m}e^{-i\o z}
{\bf f}(z,v,t)dz,\quad m=1,2. \ee
 Further, let \be\label{presU}
{\bf U}_m(\o,v,t)\defi
\E\int_t^T\exp\left\{\frac{i\o-\o^2}{2}\int_t^s\v^{v,t}(q)dq\right\}{\bf
F}_m(\o,\v^{v,t}(s),s)ds, \ee where $\v^{v,t}(s)$ is the solution
of the Ito's equation (\ref{v}) for the volatility process given
the initial condition $\v^{v,t}(t)=v$, i.e. \be\label{v2}
\left\{\ba
d_s\v^{v,t}(s )=\w\s\v^{v,t}(s) d\w w(s),\quad s>t,\\
\v^{v,t}(t)=v. \ea\right. \ee (Remind that  $\w a=0$ in (\ref{v})
by assumptions). Let \be\label{uU} {\bf
u}(z,v,t)\defi\frac{1}{\sqrt{2\pi}}\int_\R e^{i\o z}{\bf
U}_m(\o,v,t)dz, \quad z\in I_m.
\ee
\begin{theorem} \label{ThExist}  Let $K\in {\cal K}$, $\ww K=e^{-rT}K$ be given,
and let $H(\cdot)=H_j(\cdot)$ be defined as \be\label{solution}
H(x,v,t,\ww K)\defi H_{\BS,j}(x,v,t,\ww K)+{\bf
u}\left(\log\frac{x}{\ww K},v,t\right).
  \ee
Then $H(\cdot,\ww K)\in\Y(\ww K)$ and it is the unique solution of
 (\ref{parab1}) for $\varphi\equiv\varphi_j$, $j=1,2$ at this class;
 the integrals in (\ref{uU}) are
defined as the Fourier transforms of functions that are square
integrable as functions of $\o$.
\end{theorem}
{\it Remark.} The function ${\bf u}(\cdot)$ does not depend on
$j=1,2$; i.e., the correction for the Black-Scholes formula
generated by this model is the same for call and put options.
\section{Proofs}
\par
{\it Proof of Proposition \ref{propBS}}.  Consider options with
different strike prices $K\in {\cal K}$ as risky assets with
prices $S_{K,j}(t)=P_j(t,K)$.  Set $\ww S_{K,j}(t)\defi
e^{-rt}S_{K,j}(t)$. We have $P_j(t,K)=e^{rt}\ww P_j(t,\ww K)\defi
\ww H_{\BS,j}(\ww S(t),v(t),t,\ww K)$. Then
$$
\ba d\ww S_{K,j}(t)+ \frac{\p \ww H_{\BS,j}}{\p x}(\ww
S(t),v(t),t,\ww K)d\ww S(t)+ \frac{\p \ww H_{\BS,j}}{\p v}(\ww
S(t),v(t),t,\ww K)dv(t)\\+
\biggl(\frac{\p \ww H_{\BS,j}}{\p t}(\ww S(t),v(t),t,\ww K) \\
+\frac{1}{2} x^2v\frac{\p^2 \ww H_{\BS,j}}{\p x^2}(\ww
S(t),v(t),t,\ww K)+\frac{1}{2} \w\s^2v^2\frac{\p^2 \ww
H_{\BS,j}}{\p v^2}(\ww S(t),v(t),t,\ww K)\biggr)dt
\\=
\frac{\p \ww H_{\BS,j}}{\p x}(\ww S(t),v(t),t,\ww K)d\ww S(t)\\+
\frac{\p \ww H_{\BS,j}}{\p v}(\ww S(t),v(t),t,\ww
K)dv(t)+\frac{1}{2} \w\s^2v^2\frac{\p^2 \ww H_{\BS,j}}{\p v^2}(\ww
S(t),v(t),t,\ww K)dt. \ea
$$
\par Let $N\ge 2$. Let $K_i\in{\cal K}$, $\ww K_i=e^{-rT}K_i$, $i=1,2$, $K_1\neq K_2$, let $j\in\{1,2\}$
be given,   and let $X(t)$ be the wealth defined as
$$
X(t)=\b(t)B(t)+\g_0(t)S(t)+\sum_{i=1}^2\g_{i,j}(t)P_j(t,K_i).
$$
for an admissible self-financing strategy $(\b(t),\g(t))$ with
$\g(t)= (\g_0(t),\g_{1,j}(t),\g_{2,j}(t))$, such that
$\g(\cdot)\neq 0$ and
$$ \ba \g_0(t)= -\g_{1,j}(t)\frac{\p \ww H_{\BS,j}}{\p x}(\ww
S(t),v(t),t,\ww K_1)-\g_{2,j}(t) \frac{\p \ww H_{\BS,j}}{\p x}(\ww
S(t),v(t),t,\ww
K_2),\\
 \g_{1,j}(t)\frac{\p \ww H_{\BS,j}}{\p v}(\ww S(t),v(t),t,\ww
K_1)+\g_{2,j}(t) \frac{\p \ww H_{\BS,j}}{\p v}(\ww S(t),v(t),t,\ww
K_2)=0. \ea
$$ (Clearly, such a strategy exists). Then
$$
d\ww X(t)=\g_0(t)d\ww S(t)+\sum_{i=1}^2\g_{i,j}(t)d\ww
S_{K_i}(t)=\xi(t)dt,
$$
where
$$
\xi(t)=\frac{1}{2}  \g_{1,j}(t)\w\s^2v^2\frac{\p^2 \ww
H_{\BS,j}}{\p v^2}(\ww S(t),v(t),t,\ww K_1)+\frac{1}{2}
\g_{2,j}(t) \w\s^2v^2\frac{\p^2 \ww H_{\BS,j}}{\p v^2}(\ww
S(t),v(t),t,\ww K_2).
$$
We have that $\P(\xi(t)\neq 0)>0$. That means that arbitrage
possibility does exist for $N\ge 2$.\par Let $N=1$, and let
$K\in{\cal K}$, $\ww K=e^{-rT}K$.  Let $X(t)$ be the wealth
defined as
$$
X(t)=\b(t)B(t)+\g_0(t)S(t)+\sum_{j=1}^2\g_{j}(t)P_j(t,K).
$$
for an admissible self-financing strategy $(\b(t),\g(t))$ with
$\g(t)= (\g_0(t),\g_{j}(t),\g_{j}(t))$. By (\ref{BS}), it follows
that
$$\ba\frac{\p \ww H_{\BS,1}}{\p v}(\ww S(t),v(t),t,\ww
K)\equiv\frac{\p \ww H_{\BS,2}}{\p v}(\ww S(t),v(t),t,\ww
K),\\
\frac{\p^2 \ww H_{\BS,1}}{\p v^2}(\ww S(t),v(t),t,\ww
K)\equiv\frac{\p^2 \ww H_{\BS,2}}{\p v^2}(\ww S(t),v(t),t,\ww K).
\ea$$ Then
 there exist random $\F_t$-adapted processes $a_i(t)$ from $L^2(\O\times[0,T])$
such that
$$
d\ww X(t)= [\a_0(t)\g_0(t)+\a_1(t)\g_1(t)+\a_2\g_2(t)]d\ww
S(t)+[\g_1(t)+\g_2(t)[\a_2d v(t)+\a_4dt].
$$
It follows from independency of $w(\cdot)$ and $\w w(\cdot)$ that
the market is arbitrage free for $N=1$. $\Box$
 \par {\it Proof of Theorem
\ref{ThNAR}}. Let  $K\in{\cal K}$, $\ww K=e^{-rT}K$, $S_{K,j}(t)=P_j(t,K)$, $\ww
S_{K,j}(t)\defi e^{-rt}S_{K,j}(t)$ again.  We have $P_j(t,K)=e^{rt}\ww P_j(t,\ww
K)\defi \ww H_j(\ww S(t),v(t),t,\ww K)$. Then
$$
\ba d\ww S_{K,j}(t)+ \frac{\p \ww H_j}{\p x}(\ww S(t),v(t),t,\ww
K)d\ww S(t)+ \frac{\p \ww H_j}{\p v}(\ww S(t),v(t),t,\ww
K)dv(t)\\+ \left(\frac{\p \ww H_j}{\p t}(\ww S(t),v(t),t,\ww K)
+\frac{1}{2} x^2v\frac{\p^2 \ww H_j}{\p x^2}(\ww S(t),v(t),t,\ww
K)+\frac{1}{2} \w\s^2v^2\frac{\p^2 \ww H_j}{\p v^2}(\ww
S(t),v(t),t,\ww K)\right)dt
\\=
\frac{\p \ww H_j}{\p x}(\ww S(t),v(t),t,\ww K)d\ww S(t)+ \frac{\p
\ww H_j}{\p v}(\ww S(t),v(t),t,\ww K)dv(t). \ea
$$
Let $K_i\in{\cal K}$,   and let $X(t)$ be the wealth
defined
for some admissible self-financing strategy $(\b(t),\g(t))$ as
$$
X(t)=\b(t)B(t)+\g_0(t)S(t)+\sum_{j=1}^2\sum_{i=1}^N\g_{i,j}(t)P_j(t,K_i).
$$ Then
$$
d\ww X(t)=\g_0(t)d\ww S(t)+\sum_{j=1}^2\sum_{i=1}^N\g_{i,j}(t)d\ww
S_{K_i}(t).
$$
Clearly, there exist random $\F_t$-adapted processes $c_0(t)$ and $c_1(t)$ from $L^2(\O\times[0,T])$
such that
$$
d\ww X(t)= c_0(t)d\ww S(t)+c_1(t)d v(t).
$$
It follows from independency of $w(\cdot)$ and $\w w(\cdot)$ that the market is arbitrage free.
\par
{\it Remark.} Let  $K_i\in{\cal K}$, $i=1,2$.  The market with
assets ($B(t)$, $S(t)$, $S_{K_1,1}(t)$) has two risky assets and
two Brownian motions, and this market is arbitrage free. Assume
that the corresponding function $H$ in Theorem \ref{ThNAR} is such
that $\frac{\p H}{\p v}(\ww S(t),v(t),t,\ww K_1)\neq 0$. Then the
absence of arbitrage can be illustrated as \be \label{na1} \ba
d\ww S_{K_2}(t)=\frac{\p H}{\p x}(\ww S(t),v(t),t,\ww K_2)dS(t)+
\frac{\p H}{\p v}(\ww S(t),v(t),t,\ww K_2)dv(t)\\
=\frac{\p H}{\p x}(\ww S(t),v(t),t,\ww K_2)d\ww S(t)\\
+\frac{\p H}{\p v}(\ww S(t),v(t),t,\ww K_2)\left[
\frac{\p H}{\p v}(\ww S(t),v(t),t,\ww K_1)\right]^{-1}\left[d\ww S_{K_1}(t)
-\frac{\p H}{\p v}(\ww S(t),v(t),t,\ww K_1)d\ww S(t)\right].
\ea
\ee
 By (\ref{na1}), it follows that
 for any strategy for the market with  assets ($B(t)$, $S(t)$, $S_{K_1,1}(t)$,
 $S_{K_2,1}(t)$), there exists a
strategy  using assets ($B(t)$, $S(t)$, $S_{K_1,1}(t)$) only with
the same wealth. Thus, adding options with different strike prices
does not lead to arbitrage.
\par
{\it Proof of Theorem \ref{ThExist}}. Instead of solving the
problem (\ref{parab1}), we shall solve the  problem \be
\label{parab12} \left\{ \ba \frac{\p G}{\p t}(x,v,t) +\frac{1}{2}
x^2v\frac{\p^2 G}{\p x^2}(x,v,t)+\frac{1}{2} \w\s^2v^2\frac{\p^2
G}{\p v^2}(x,v,t)=\phi(x,v,t),
\\
G(x,v,T)=0,\\
G(\ww K,v,t)=0,\\
\left.\frac{\p G }{\p x}(x,v,t)\right|_{x=\w K}=0,\ea \right. \ee
given function $\phi:Q\to\R$.
 Clearly, the desired solution  $H(x,v,t,\ww K)$ can be found as $H(x,v,t,\ww K)=\ww H_\BS
(x,v,t,\ww K)+G(x,v,t)$ if $
\phi(x,v,t)=\frac{1}{2}\w\s^2v^2\frac{\p^2\ww H_\BS}{\p
v^2}(x,v,t)$.
\par
Introduce the Banach space $\X$ of functions $\phi:Q\to\R$ such
that $$\ba \|\phi(\cdot)\|_\X&\defi\sup_{(v,t)\in Q}
v^{-1/2}\left(\int_\R|f(z,v,t)|dz+\left[\int_\R|f(z,v,t)|^2dz\right]^{1/2}\right)<+\infty,
\\& \hbox{where}\quad f(z,v,t)\defi \phi(\ww K e^z,v,t). \ea
$$
\par
  Let $j\in\{1,2\}$ and $\ww K$ be fixed. Let $D_1(\ww K)\defi
\{(x,v)\in D:\ \ x<\ww K\}$, $D_2(\ww K)\defi \{(x,v)\in D:\ \
x>\ww K\}$. Let $Q_m(\ww K)\defi D_m(\ww K)\times [0,T]$, $m=1,2$.
Let ${\cal Y}_{j,m}(\ww K)$ denotes the class of functions
$u=u(x,v,t):Q_m(K)\to\R$ such that $u(\cdot)=\ww
u(\cdot)|_{Q_m(K)}$, where $\ww u(\cdot)\in \Y_j(\ww K)$, $j=1,2$.
\par
 Instead of solving the problem (\ref{parab1}), we shall
investigate the pair of problems in $Q_m(\ww K)$, and $m=1,2$.
Clearly, it suffices to prove the existence of a solution
$G(\cdot)\in {\cal Y}_{1,m}(\w r,\ww K)$  of  the problem
(\ref{parab12}) in $Q_m(\ww K)$, $m=1,2$. (Recall that $\ww
K=e^{-rT}K$).
\par
  Let $\phi\in \X$, \be\label{oou}
  \u(z,v,t)\defi G(x,v,t),\quad  f(z,v,t)\defi \phi(x,v,t),\quad
   z\in \R, \quad x=\ww K e^z.
  \ee
  Formally,
  $$\frac{\p G}{\p x}(x,v,t)=\frac{1}{x}\frac{\p \u}{\p z}(z,v,t),\quad
\frac{\p^2 G}{\p x^2}(x,v,t)=\frac{1}{x^2}\frac{\p^2 \u}{\p
z^2}(z,v,t) -\frac{1}{x^2}\frac{\p \u}{\p z}(z,v,t).
  $$
Set $\w Q_1\defi \{(z,v,t): \ z\le 0\}$, $\w Q_2\defi \{(z,v,t): \
z\ge 0\}$. Problem (\ref{parab12}) can be rewritten for $\u:\w
Q_m\to\R$, $m=1,2$, as \be \label{parabu} \left\{ \ba \frac{\p
\u}{\p t}(z,v,t) +\frac{v}{2}\frac{\p^2 \u}{\p
z^2}(z,v,t)-\frac{v}{2}\frac{\p \u}{\p z}(z,v,t)+\frac{1}{2}
\w\s^2v^2\frac{\p^2 \u}{\p v^2}(z,v,t)= f(z,v,t),
\\
\u(z,v,T)=0\\
\u(0,v,t)=0\\
\frac{\p \u}{\p z}(0,v,t)=0. \ea \right. \ee
\par
 Assume that we can
find the solution $u(z,v,t)$ of this problem in the class of
functions  that belong $L_2(\R)$ as functions of $z$ and the
following Fourier transform is defined: \be\label{U} \ba
U(\o,v,t)\defi\frac{1}{\sqrt{2\pi}}\int_{I_m}e^{-i\o z}u(z,v,t)dz,
\quad\o\in\R, \quad m=1,2.\ea \ee
 Instead of  (\ref{parab12}), consider the problem for
$U$ \be \label{parabU} \left\{ \ba \frac{\p U}{\p t}(\o,v,t)
+\left(-\o^2-i\o \right)\frac{v}{2}U(\o,v,t)+\frac{1}{2}
\w\s^2v^2\frac{\p^2 U}{\p v^2}(\o,v,t)=F(\o,v,t),\\
U(\o,v,T)=0, \ea \right. \ee where \be\label{Ff2}
F(\o,v,t)\defi\frac{1}{\sqrt{2\pi}}\int_{I_m}e^{-i\o
z}f(z,v,t)dz.\ee
\begin{proposition}\label{pVS1}  There exist a constant $C_1>0$ such that
for any $\phi(\cdot)\in \X$
$$\ba
\left|F(\o,v,t)\right|\le C_1\sqrt{v}\|\phi(\cdot)\|_\X\quad
\forall z\in\R,\ v>0,\ t\in[0,T],\ \o\in\R, \ea
$$ where $F(\cdot)$ is defined by (\ref{Ff2}) and (\ref{oou}).
\end{proposition}
 \par
 {\it Proof.} We have that
 $f(z,v,t)=c(z,v,t)\sqrt{v},$ where
$\int_{I_m}|c(z,v,t)|dz\le \|\phi(\cdot)\|_\X$. Clearly,   the
Fourier transform of $c(\cdot,v,t)$ is uniformly bounded. $\Box$
\par Further, the solution of (\ref{parabU}) given $\o\in\R$
exists and can be expressed as \be\label{pres}
 U(\o,v,t)=
\E\int_t^T\exp\left\{\frac{-\o^2+i\o}{2}\int_t^s\v^{v,t}(q)dq\right\}F(\o,\v^{v,t}(s),s)ds,
\ee where $\v^{v,t}(s)$ is the solution of the Ito's equation
(\ref{v2}). \par
\begin{lemma} \label{lemmaest} There exist a constant $C>0$ such that
 \be\label{estU}
\ba\int_{-\infty}^{+\infty}(1+|\o|^{2})|U(\o,v,t)|^2d\o\le
C\|\phi(\cdot)\|_\X \quad \forall \ v>0,\ t\in[0,T],
\\
 \int_{-\infty}^{+\infty}(1+|\o|^2)^{p\d /2}|U(\o,v,t)|^{p}d\o\le C\|\phi(\cdot)\|_\X
\quad \forall\ p\ge 1,\ \d\in [0,2-1/p), \ v>0. \ea \ee
\end{lemma}
\par {\it Proof}. We have that
$$
\ba |U(\o,v,t)|&\le \frac{\ww K}{2}
\E\int_t^T\left|F(\o,\v^{v,t}(s),s)\right|
\exp\left\{\frac{-\o^2}{2}\int_t^s\v^{v,t}(q)dq\right\} ds
\\
&\le\const\|\phi(\cdot)\|_\X\E\int_t^T\sqrt{\v^{v,t}(s)}\exp\left\{\frac{-\o^2}{2}
\int_t^s\sqrt{\v^{v,t}(q)}dq\right\}ds\\
&\le\const\|\phi(\cdot)\|_\X\sqrt{T}\,\E\int_t^T\v^{v,t}(s)\exp\left\{\frac{-\o^2+\a_0}{2}
\int_t^s\v^{v,t}(q)dq\right\}ds\\&=
\const\|\phi(\cdot)\|_\X\frac{2}{\o^2}\E\left(1-\exp\left\{\frac{-\o^2}{2}
\int_t^T\v^{v,t}(q)dq\right\}\right)
 \ea
$$
for all $v>0,\ t\in[0,T]$. This completes the proof. $\Box$
\par
 Let $$\rho(y)\defi\frac{1}{1+e^{3y}},\quad  y\in\R.
$$
\par Consider the Hilbert space
$L_2(\R,\rho)$ with the weight $\rho(\cdot)$ that consist of
complex valued functions $\eta:\R\to\R$ such that
$\|\eta(\cdot)\|_{L_2(\R,\rho)}\defi
\left(\int_\R|\eta(y)|^2\rho(y)dy\right)^{1/2}<+\infty$. The
scalar product for that spaces is
$$
(\eta_1(\cdot),\eta_2(\cdot))_{L_2(\R,\rho)}\defi \int_\R \oo
\eta_1(y)\eta_2(y)\rho(y)dy.
$$
Let $\oo X$ be the Banach space of functions $\xi(\cdot):\R\times
[0,T]\to\R$  such that $\xi(\cdot)\in C([0,T];L_2(\R,\rho))$ and
$$
\|\xi(\cdot)\|_{\oo
X}\defi\sup_{t\in[0,T]}\|\xi(\cdot,t)\|_{L_2(\R,\rho)} +
\left(\int_0^T\left\|\frac{\p \xi}{\p y}(\cdot,t)
\right\|_{L_2(\R,\rho)}^2dt\right)^{1/2}<+\infty. $$\par
 Set
\be\label{psi}
\ba \psi(z,y,t)\defi u(z,e^y,t)=G(\ww K e^z,e^y,t),\\
\Psi(\o,y,t)\defi U(\o,e^y,t),\\
 h(\o,y,t)\defi F(\o,e^y,t).
\ea
\ee
\begin{lemma} \label{lemmaestV}
There exist a constant $C>0$ such that \be\label{estU3} \ba
\int_\R\|V(\o,\cdot)\|^2_{\oo X}d\o\le C \|\phi(\cdot)\|_\X.\ea
\ee
\end{lemma}
\par
{\it Proof.} Set
$$
A\Psi=A(\o,y,t)\Psi\defi \frac{\w\s^2}{2} \frac{\p^2 \Psi}{\p
y^2}(\o,y,t)-\frac{\w\s^2}{2}\frac{\p \Psi}{\p y}(\o,y,t), \quad
g(\o,y,t)=-\frac{e^y\o^2}{2}-\frac{e^yi\o}{2}.
$$ By the definitions, the function $\Psi$ satisfies \be \label{parabPsi} \left\{ \ba
\frac{\p \Psi}{\p t}(\o,v,t)
+A\Psi(\o,y,t)+g(\o,y,t)\Psi(\o,y,t)=h(\o,y,t),\\
\Psi(\o,y,T)=0. \ea \right. \ee Let the Hilbert space $\L$ be
defined as $\L\defi L^2([0,T];L_2(\R,\rho))$.
\par
Let $\{h^{(j)}(\o,y,t)\}_{j=1}^{+\infty}$ be a sequence of bounded
smooth functions such that $h^{(j)}\to h$ in
$\L=L^2([0,T];L_2(\R,\rho))$ for any $\o$. Let
$\Psi=\Psi^{(j)}(\o,\cdot)$ be the classical solution of
(\ref{parabU}) given $\o$ for $h=h^{(j)}$.  Clearly,
$$
\Re\bigl(\Psi(\o,\cdot,t),
g(\o,\cdot,t)\Psi(\o,\cdot,t)\bigr)_{L_2(\R,\rho)}=\Re\int_\R|\Psi(\o,y,t)|^2
\frac{e^y(-\o^2-i\o)}{2}\rho(y)dy\le 0.$$ Hence \be\label{big}\ba
&\frac{1}{2}\|\Psi(\o,\cdot,s)\|^2_{ L_2(\R,\rho)}-\frac{1}{2}\|
\Psi(\o,\cdot,T)\|^2_{L_2(\R,\rho)}\\&=\Re\int^{T}_s
(\Psi(\o,\cdot,t),A\Psi(\o,\cdot,t)
+g(\o,\cdot,t)\Psi(\o,\cdot,t)+h(\o,\cdot,t) )_{ L_2(\R,\rho)} dt
\\&=\Re\int^{T}_s \biggl\{-\frac{\w\s^2}{2}\biggl( \frac{\p
\Psi}{\p y}(\o,\cdot,t),\frac{\p \Psi}{\p y}(\o,\cdot,t)
\biggr)_{L_2(\R,\rho)} \\& -\biggl( \Psi(\o,\cdot,t),\frac{\p
\rho}{\p y}(\cdot)\frac{1}{\rho(\cdot)} \frac{\p \Psi}{\p
y}(\o,\cdot,t) \biggr)_{L_2(\R,\rho)}+\biggl(
\Psi(\o,\cdot,t),\frac{\p \Psi}{\p y}(\o,\cdot,t)
\biggr)_{L_2(\R,\rho)} \\
&\hphantom{xxxxxxxxxxxx}+\bigl(\Psi(\o,\cdot,t),
g(\o,\cdot,t)\Psi(\o,\cdot,t)+h(\cdot,t)
\bigr)_{L_2(\R,\rho) }\biggr\} dt\\
&\le \int^{T}_s    \biggl\{
-\left(\frac{\w\d^2}{2}-\e\right)\biggl\|\frac{\p \Psi}{\p y}
(\o,\cdot,t)\biggr\|^2_{L_2(\R,\rho)}\\
&\hphantom{xxxx}+c_1(\e)\| \Psi(\o,\cdot,t)\|^2_{ L_2(\R,\rho)}
+\|\Psi(\o,\cdot,t)\|_{L_2(\R,\rho)}\|h(\o,\cdot,t)\|_{L_2(\R,\rho)}\biggr\}
dt \\ &\le \int^{T}_s \biggl\{ -\left(\frac{\w\d^2}{2}-\e\right )
\biggl\|\frac{\p \Psi}{\p y} (\o,\cdot,t)\biggr\|^2_{L_2(\R,\rho)}
\\&\hphantom{xxxxxxxxxxxxxxx}+c_2(\e)\|
\Psi(\o,\cdot,t)\|^2_{ L_2(\R,\rho)}
+c_3(\e)\|h(\o,\cdot,t)\|^2_{L_2(\R,\rho)}\biggr\} dt. \ea \ee
Here $\e \in (0,\w\d^2/2 )$ can be any, and  $c_k(\e )$ depends
only on $\e$; we have used the inequality $2\alpha \beta \le \e
\alpha +\e^{-1}\beta$ $(\forall \alpha ,\beta ,\e \in \R, \e
>0)$. To derive (\ref{big}), we have also used  the equality
$$ \biggl( \Psi(\o,\cdot,t),\frac{\p \rho}{\p y}(\cdot) \frac{\p \Psi}{\p
y}(\o,\cdot,t) \biggr)_{L_2(\R)} =\left( \Psi(\o,y,t),\frac{\p
\rho}{\p y}(y)\frac{1}{\rho(y)} \frac{\p \Psi}{\p y}(\o,y,t)
\right)_{L_2(\R,\rho)}.
 $$
Note that the function $\frac{\p \rho}{\p y}(y)\frac{1}{\rho(y)}$
is bounded. \par By Bellman inequality, it follows from
(\ref{big}) that there exists a constant $C>0$ such that
 \be\label{estU33}
\ba\sup_{t\in[0,T]}\|\Psi(\o,\cdot,t)\|^2_{L_2(\R,\rho)} +
\int_0^T\left\|\frac{\p \Psi}{\p y}(\o,\cdot,t)
\right\|^2_{L_2(\R,\rho)}dt\le C\|h(\o,\cdot)\|^2_\L \quad \forall
\o\in\R. \ea \ee Hence there exists a constant $C>0$ such that
 \be\label{estU34}
\|\Psi(\o,\cdot,t)\|_{\oo X}\le C\|h(\o,\cdot)\|_\L \quad \forall
\o\in\R. \ee Remind that $\sup_{\o,t}\|h(\o,\cdot)\|_\L\le \const
\|\phi\|_\X$, and $\Psi=\Psi^{(j)}$ is the solution of
(\ref{parabU}) for $h=h^{(j)}$. Clearly, the sequence
$\{\Psi^{(j)}\}_{j=1}^{+\infty}$ has the limit $\Psi\in\oo X$ that
we shall call the solution of (\ref{parabU}) for the original $h$
(this is in fact a generalized solution; the uniqueness of this
solution follows from the linearity of the problem and from
(\ref{estU34})). \par Further, it follows from the estimation
$$\int_\R|F(\o,v,t)|^2dw=\const\int_\R|f(z,v,t)|^2dz\le
\const v\|\phi(\cdot)\|^2_\X$$ that
$$\int_\R\|h(\o,\cdot)\|^2_{\L}d\o=\int_\R d\o \int_0^T
\|h(\o,\cdot,t)\|^2_{L_2(\R,\rho)}dt\le\const\|\phi(\cdot)\|^2_\X.
$$ This completes the proof of Lemma \ref{estU3}. $\Box$
\par
Remind that $\Psi(\o,y,t)=U(\o,e^y,t)$. Clearly, $$ \ba \frac{\p
\Psi}{\p y}(\o,y,t)=\frac{\p U}{\p y}(\o,e^y,t)e^{y},\\
\psi(z,y,t)=\frac{1}{\sqrt{2\pi}}\int_\R e^{i\o z}\Psi(z,y,t)dz,
\\
\frac{\p \psi}{\p y}(z,y,t)=\frac{1}{\sqrt{2\pi}}\int_\R e^{i\o
z}\frac{\p \Psi}{\p y}(z,y,t)dz. \ea
$$
By Lemma
 \ref{estU3},
 \be\label{estU4}
\int_0^Tdt\int_{I_m}dz\int_{\R}\rho(y)\left(|u(z,e^y,t)|+\left|\frac{\p
u}{\p y}(z,e^y,t)e^{y}\right|\right)^2 dz
\le\const\|\phi(\cdot)\|^2. \ee Further, it follows from Lemma
\ref{lemmaest}  that
 \be\label{estU4cor}
\sup_{t\in[0,T]}\int_{\R}
e^{-|y|}dy\int_{I_m}\left(|u(z,e^y,t)|+\left|\frac{\p u}{\p
z}(z,e^y,t)\right|\right)^2 dz \le\const\|\phi(\cdot)\|^2_\X. \ee
\begin{lemma}\label{cont} For any $\phi\in\X$, the derivative $\p u(z,v,t)/\p z$
is continuous in $z$ for a.e. $v,t$, and $\esssup_{y,t}\left|
\frac{\p \psi}{\p z}(z,y,t)\right|\le\const\|\phi(\cdot)\|_\X$.
\end{lemma}
\par
{\it Proof}. Remind that our selection of the version of the
Fourier transforms means that  $u(\cdot,v,t)$ is extended by zero
from $I_m$ (or $u(\cdot)$ is extended by zero from $\w Q_m$ to $\w
Q_{3-m}$, $m=1,2$, where $\w Q_1\defi \{(z,v,t): \ z\le 0\}$, $\w
Q_2\defi \{(z,v,t): \ z\ge 0\}$).  By Lemma \ref{lemmaest}, it
follows that $u(\cdot,v,t)\in L^{p,\d}(\R)$ and
$\esssup_{y,t}\left\|\frac{\p \psi}{\p
z}(\cdot,y,t)\right\|_{L^{p,\d}(\R)}\le\const\|\phi(\cdot)\|_\X$
with $p=3$ and $\d\in (\frac{4}{3},\frac{5}{3})$, where
$L^{p,s}(\R)$ is the space described in Adams (1975), p. 220: it
is the image of $L^p(\R)$ under the linear mapping $J^\d$ such
that $J^su={\bf F}^{-1}(1+|\cdot|^2)^{-\d/2}{\bf F}u)$, where
${\bf F}$ is the Fourier transform. By imbedding Theorems 7.63(g)
and 7.57(c) from Adams (1975), it follows that the derivatives $\p
u(z,v,t)/\p z$ are continuous for a.e. $v,t$ and
$$
\esssup_{y,t}\left(|\psi(z,y,t)|+\left|\frac{\p \psi}{\p
z}(z,y,t)\right|\right)\le\const\|\phi(\cdot)\|_\X . $$ $\Box$
\begin{lemma}\label{cont2} For any $\phi\in\X$, the function $ u(z,v,t)$ is
continuous in $t$ for a.e. $z,t$.
\end{lemma}
\par
{\it Proof}. It follows from (\ref{pres}) that $U(\o,v,t)$ is
continuous in $t$ for all $\o,v$. Further, it follows from Lemma
\ref{lemmaest} that $\|U(\cdot,v,t)\|_{L_1(\R)}\le
\const<+\infty$. By Dominated Convergency Theorem it follows
continuity of the Fourier transform.  $\Box$
\par
We can summarize the results of Lemmas \ref{lemmaest}-\ref{cont2}
in the following theorem.
\begin{theorem}\label{thest} For any $\phi\in\X$, the function $\psi(z,y,t)$, defined in (\ref{psi})
belongs $\Y$,
and there exists a constant $C>0$ such that \be\label{energ}
\|\psi(\cdot)\|_\Y\le\const\|\phi(\cdot)\|_\X \quad \forall \phi\in
\X. \ee
\end{theorem}
\begin{proposition}\label{properBS} In (\ref{varphi}),
$\varphi(x,v,t)=\frac{1}{2}\w\s^2v^2\frac{\p^2\ww H_\BS}{\p
v^2}(x,v,t),$ and $\varphi(\cdot)\in \X$.
\end{proposition}
\par
{\it Proof.} By (\ref{BS}), it follows that
$$
\ba \frac{\p\ww H_{\BS}}{\p v}(x,v,t,\ww K)&= \frac{x
}{\sqrt{2\pi}} e^{-\frac{d_+(x,v,t,\ww K)^2}{2}}\frac{\p d_+}{\p
v}(x,v,t,\ww K)- \frac{\ww K }{\sqrt{2\pi}}
e^{-\frac{d_-(x,v,t,\ww K)^2}{2}}\frac{\p\oo d}{\p v}(v,t)\\
&=\frac{x}{\sqrt{2\pi}}e^{-\frac{d_+(x,v,t,\ww K
)^2}{2}}\left[-\frac{\log{x}-\log{\ww
K}}{2v\sqrt{(T-t)v}}+\frac{\sqrt{(T-t)}}{4\sqrt{v}}\right]\\
&\hphantom{xxxxxxx}-\frac{\ww K}{\sqrt{2\pi}}e^{-\frac{
d_-(x,v,t,\ww K )^2}{2}}\left[-\frac{\log{x}-\log{\ww
K}}{2v\sqrt{(T-t)v}}-\frac{\sqrt{(T-t)}}{4\sqrt{v}}\right]. \ea
$$
Then
$$ \ba \frac{\p^2\ww H_{\BS}}{\p v^2}(x,v,t,\ww K)\\
=\frac{x}{\sqrt{2\pi}}e^{-\frac{d_+(x,v,t,\ww K
)^2}{2}}\biggl(\left[-\frac{\log{x}-\log{\ww
K}}{2v\sqrt{(T-t)v}}+\frac{\sqrt{(T-t)}}{4\sqrt{v}}\right]^2
+\frac{3}{2}\frac{\log{x}-\log{\ww
K}}{2v^2\sqrt{(T-t)v}}-\frac{\sqrt{(T-t)}}{8v\sqrt{v}}
\biggr)\\
\hphantom{xxx}-\frac{\ww K}{\sqrt{2\pi}}e^{-\frac{ d_-(x,v,t,\ww K
)^2}{2}}\biggl(\left[-\frac{\log{x}-\log{\ww
K}}{2v\sqrt{(T-t)v}}-\frac{\sqrt{(T-t)}}{4\sqrt{v}}\right]^2+
\frac{3}{2}\frac{\log{x}-\log{\ww
K}}{2v^2\sqrt{(T-t)v}}+\frac{\sqrt{(T-t)}}{8v\sqrt{v}}\biggl). \ea
$$
Clearly, $\varphi(\cdot)$ is continuous in $(0,+\infty)\times
(0,+\infty)\times [0,T)$, and the limit of $\varphi(\cdot)$ is
${\cal O}(v^{1/2})$  as $x\to 0$, $x\to \ww K$, $x\to +\infty$ as
well as  $v\to 0$, $v\to +\infty$ and  $t\to T$. The required
summarability is ensured by properties of $e^{d_\pm}$. $\Box$
\par Now we are in the position to complete the proof of Theorem
\ref{ThExist}.
\par From Theorem
\ref{thest}, the existence of solution in Theorem \ref{ThExist}
follows. The fact that the solution is unique in the given class
follows from (\ref{energ}). The formula for ${\bf u}$ from
(\ref{solution}) is in fact the formula derived above for $\u$
with substituting $\phi(\cdot)=\varphi( \cdot)$. $\Box$
\section*{ References} \ \par Adams, R.A. (1975).
{\it Sobolev Spaces.} Academic press.
\par
Black, F. and M. Scholes (1972): The valuation of options contracts
and test of market efficiency. {\it Journal of Finance}, {27},
399-417.
\par
Black, F. and M. Scholes (1973): The pricing of options and
corporate liabilities. {\it Journal of Political Economics}, {81}, 637-659.
\par Christie, A. (1982): The stochastic behaviour of
common stocks variances: values, leverage, and interest rate
effects. {\it Journal of Financial Economics}, {10},
407-432.
\par Day, T.E. and C.M. Levis (1992): Stock
  market volatility and the information
content of stock index options. {\it Journal of Econometrics}, {52}, 267-287.
\par
Derman, E., I. Kani,  and J.Z. Zou (1996): The local volatility
surface: unlocking the information in index option prices. {\it
Financial Analysts Journal} 25-36.\par
Geman, H. and T. Ane (1996):
Stochastic subordination. {\it Risk}, {9} (9), 145-149.
 \par
   Dokuchaev, N.G. (1995):  Probability
distributions  of  Ito's processes: estimations for density
functions and for conditional expectations of integral
functionals. {\it Theory of Probability and Its
Applications,} 39, iss. 4, 662-670.
\par
Dokuchaev, N.G., and X.Y. Zhou (2001): Optimal investment strategies with bounded
risks, general utilities, and goal achieving. {\it Journal of Mathematical Economics},
35, iss.2,   289-309.\par
 Dokuchaev N.G. (2002): {\it Dynamic portfolio strategies:
quantitative methods and empirical rules for incomplete information.}
 Kluwer Academic Publishers, Boston.
\par
Harrison, J.M., and S.R. Pliska (1981): Martingales and stochastic integrals in
the theory of stochastic trading. {\it Stochastic Processes and their
Application}, { 11}, 215-260.
\par Hauser, S. and B. Lauterbach (1997): The relative performance of five
alternative warrant pricing models. {\it Financial Analysts Journal},  N1, 55-61.
\par
Hull, J. and A. White (1987): The pricing of options on assets with
stochastic volatilities. {\it Journal of Finance}, { 42},
281-300.\par Jarrow, R. (ed.) (1998): {\it Volatility New Estimations Techniques for
pricing Derivatives}, Risk Books.
\par Johnson, H. and D. Shanno
(1987): Option pricing when the variance is changing. {\it Journal of Financial and
Quantitative Analysis}, { 22}, 
143-151.
\par
Jouini, E. Market imperfection, equilibrium and
arbitrage. (1996):
{\it Financial Mathematics} ({\it Lecture Notes in Mathematics},
{ 1656}), 247-307.
\par
  Masi, G.B., Yu.M. Kabanov,  and W.J. Runggaldier (1994):
 Mean-variance hedging of options on stocks
with Markov volatilities. {\it Theory of Probability and Its
Applications
} { 39}, 
172-182.
\par Mayhew, S. (1995): Implied volatility. {\it Financial
Analysts Journal}, iss. 4,  8-20.
\par
 Taylor, S.J. and  X. Xu (1994): The magnitude of
implied volatility smiles: theory and empirical evidence for
exchange rates.
{\it Review of Future Markets}, { 13}, 
355-380.
\end{document}